\documentclass[journal]{IEEEtran}

\usepackage{graphicx} 
\usepackage{amsmath}
\usepackage{amsthm,amsfonts}
\usepackage{amssymb}
\usepackage[greek,english]{babel}
\usepackage{enumitem}
\usepackage{cite}

\usepackage{pgfplots}
\usepackage{tikz}
\usepackage[customcolors]{hf-tikz}
\usepgfplotslibrary{colorbrewer}
\usetikzlibrary{3d}

\newcommand{\gpi}{\textrm{\greektext p}}
\renewcommand{\pi}{\gpi}
\providecommand*{\M}[1]{\mathbf#1} 
\providecommand*{\mrm}[1]{\mathrm{#1}}

\providecommand*{\V}[1]{\boldsymbol#1} 
\renewcommand{\vec}[1]{\V{#1}} 
\providecommand*{\UV}[1]{\hat{\boldsymbol#1}}
\providecommand*{\T}[1]{\mathrm{#1}} 

\providecommand*{\eu}{\ensuremath{\mrm{e}}}

\providecommand*{\ju}{\ensuremath{\mrm{j}}}
\providecommand*{\diff}{\operatorname{d}\!}

\tikzset
{%
  front face/.style={fill=green!50,canvas is xy plane at z=-\Rx},
  up    face/.style={fill=green!50,canvas is xz plane at y=\h},
  east  face/.style={fill=green!50,canvas is yz plane at x=-\Rx},
  pics/square/.style={
    code={\draw[fill=green!30,even odd rule] (-\Rt,-\Rt) rectangle (\Rt,\Rt) (-\Rx,-\Rx) rectangle (\Rx,\Rx);}},
  pics/sida/.style={
    code={\draw[fill=green!30,even odd rule] (-\Rt,0) rectangle (\Rt,\h);}},		
  pics/sidas/.style={
    code={\draw[fill=green!30,even odd rule] (0,-\Rt) rectangle (\h,\Rt);}},		
}
\tikzset{>=latex}

\title{Interpreting Moment Matrix Blocks Spectra using Mutual Shadow Area}

\author{
Yaniv~Brick,
Francesco P.~Andriulli,
and Mats~Gustafsson
\thanks{Manuscript received \today; revised \today. This work was supported in part by the Israel Science Foundation (ISF) grant No. 442/22, the TICRA foundation, the Swedish Research Council SEE-6GIA, and SSF Sabbatical.}
\thanks{Y. Brick is with the School of Electrical and Computer Engineering, Ben-Gurion University of the Negev, Beer-Sheva, Israel 8410501 (e-mail:
ybrick@bgu.ac.il). 

F. P. Andriulli is with the Department of Electronics and Telecommunication, Politecnico di Torino, Turin, Italy, 10129 (e-mail:
francesco.andriulli@polito.it). 

M. Gustafsson is with Electrical and Information Technology, Lund University, Lund, Sweden, (e-mails: mats.gustafsson@eit.lth.se).}
}

\begin{document}

\maketitle
\begin{abstract}
The mutual shadow area of pairs of surface regions is used for guiding the study of the spectral components and rank of their wave interaction, as captured by the corresponding moment matrix blocks. 
It is demonstrated that the mutual shadow area provides an asymptotically accurate predictor of the location of the singular value curve knee. 
This predicted knee index is shown to partition the interacting parts of the range and domain of blocks into two subspaces that can be associated with different wave phenomena: an ``aperture" subspace of dimension that scales with the subdomains area (or length in 2-D) and a remainder ``diffraction" subspace of dimension that scales much slower with the electrical length, depending on the geometric configuration.  
For interactions between open surface domains typical for the common hierarchical partitioning in most fast solvers, the latter can be attributed to the domain edges visible by its interacting counterpart. 
For interactions in 3-D with a small aspect angles between the source and observers, the diffraction subspace dimension is dominant in determining the rank until fairly large electrical lengths are reached. This explains
the delayed asymptotic scaling of ranks and impressive fast solver performance observed in recent literature for seemingly arbitrary scatterers with no special geometric characteristics. 
In the extreme cases of ``endfire" reduced dimensionality interactions, where the shadow area vanishes, the diffraction governs also the asymptotic rank, which translates to superior asymptotic solver performance. 
\end{abstract}


\section{Introduction}
The concept of degrees of freedom (DoF) of wave interactions is of importance in various contexts. 
In the analysis of radiation and scattering, it is related to the number of sampling points, fundamental wave types, and characteristic modes required for accurately capturing the radiated or scattered near and far fields~\cite{Bucci+Isernia1997,Franceschetti2017,Bucci+Migliore2025,Gustafsson+Lundgren2024,DiFrancia1955,Maisto+etal2021,Miller2019}. 
It is also related to the concept of channel capacity and effective number of communication and radar channels, specifically in multiple-input multiple-output systems, both in free space~\cite{Pizzo+etal2022,Poon+etal2005,Dardari2020,Puggelli+etal2025} and within wave-guiding structures~\cite{Arendt+etal2009}. 

In computational electromagnetics and acoustics, the number of DoF has been discussed in the context of fast integral equation (IE) solvers. 
Such solvers make use of efficient representations of wave interactions, through field integral operators, between subsets of basis and testing functions confined to finite-sized and arbitrary-shaped subdomains of an analyzed system ~\cite{Chew+Lu1993,Bebendorf2000,hackbuschSparseMatrixArithmetic1999,bormDatasparseApproximationNonlocal2007}.
The number of DoF is thought of as the number of different wave types sufficient for describing the interactions. 
As such, it is related (but not identical) to those of multipoles or plane waves required for the accurate representation of the fields radiated by a finite-support source distribution~\cite{Chew1995,songMultilevelFastmultipoleAlgorithm1995,coulierInverseFastMultipole2017} or to that of sampling points from which they can be reconstructed through interpolation in a given domain \cite{bucciDegreesFreedomScattered1989, boagNonuniformPolarGrid2002, brickMultilevelNonuniformGrid2010, brickFastDirectSolution2011}. 
The observation that, for domains interacting in 3-D, these numbers scale asymptotically as \((ka)^2\) (where \(k\) is the wavenumber and \(a\) is a typical dimension of the  domains, such that \(ka\) loosely denotes their electrical lengths) has driven the development of ``matrix-free" fast solvers, in which various economical descriptions of partial contributions to the interactions are used to hierarchically represent the moment matrix.  

Moment matrix blocks can be also expressed, to a prescribed accuracy, by their algebraically computed low-rank (LR) approximations.  
Like the number of DoF, the rank required for guaranteeing a prescribed relative error level for a block is dictated by the geometric properties of the participating subdomains. 
 Numerous fast solvers rely on the partitioning of the moment matrix into hierarchies of blocks, of which most are algebraically compressed using LR approximation \cite{hackbuschSparseMatrixArithmetic1999, bormDatasparseApproximationNonlocal2007, zhaoAdaptiveCrossApproximation2005, martinssonFastDirectSolver2005,Greengard2012, weiFastDirectMatrix2012, ambikasaranMathcalLogFast2013,  laiFastDirectSolver2014, brickFastDirectSolver2014, bormDirectionalmatrixCompressionHigh2017}. 
 The costs of building, storing, applying to vectors, and factorizing these representations and their asymptotic scaling with the number of problem unknowns \(N\) depend greatly on the ranks. 
 For objects and error thresholds where the largest block ranks scale slower than \(N\) with \(ka\), asymptotic savings can be achieved. 
With multipoles and spatial bandwidth analyses providing only loose bounds on the {\it{rank}}, the scaling of ranks has been studied numerically.
This has been done explicitly, as in \cite{brickFastMultilevelComputation2016, brickRapidRankEstimation2018, kelleyIterativeRandomSampling2023}, as well as implicitly through the performance of various algebraic compression schemes and fast solvers. 

 For ``endfire" interactions, which are of reduced dimensionality (e.g., in elongated or quasi-planar settings), slow asymptotic scaling with \(ka\) has been demonstrated repeatedly \cite{michielssenScatteringElongatedObjects1996, martinssonFastDirectSolver2005, martinssonFastDirectSolver2007}, including for specialized modified IE kernels designed for reducing the effective dimensionality of interactions \cite{brick_fast_2014a, sharshevskyDirectSolutionScattering2020,zvulunGeneralizedSourceIntegral2023 ,dahanFastDirectSolvers2024, dahanVectorGeneralizedSource2025}.
 However, conventional kernel interactions between regions of arbitrary shape objects may include broadside components. 
 For these, the scaling may appear surprisingly slow with \(N\) until large \(ka\) values are reached  \cite{brickFastMultilevelComputation2016,brickRapidRankEstimation2018}. 
 In some geometric and threshold configurations, an eventual \(\mathcal{O}(N)\) scaling is observed. 
 In others, particularly for low threshold values and relatively remote interacting subdomains, deterioration of the scaling shows only for subdomains of at least dozens of wavelengths in size, contradicting the intuition gained from spectral and spatial sampling bounds.
 As reasonably sized realistic complex objects involve mixtures of broadside and endfire interactions between domains of various sizes, the extraction of 
 asymptotic trends from compression and performance data for the corresponding entire moment matrices can be misleading.
 The impressive performance of fast LR compression-based direct solvers for such objects, combined with the insidious behavior of the rank, loose analytical bounds on it, and an ambiguous definition of the number of DoF, has sustained debates on the expected high-frequency asymptotic scaling of ranks and solver performance. 

Recently, an approach for estimating the asymptotic number of DoF of the far fields radiated by electrically large source domains was presented~\cite{Gustafsson2025a}. 
It has been shown that the knee index in the singular value (SV) plot describing the radiation by a domain is predicted accurately by using the cumulative shadow area (or, in 2-D, length) cast by the object when illuminated by all possible plane waves. 
For the interaction between finite-size source and observation regions, this concept can be extended to a cumulative {\it{mutual}} shadow area~\cite{Gustafsson2025c}. 
This value can be calculated quickly, at an \(\mathcal{O}(1)\) cost, based solely on geometric considerations. 
For some special cases, it even has a closed-form expression. 
The prediction of the knee index is highly relevant for communication channel capacity assessment, where it represents well the dimension (or number of DoF) of the signal space. 
In the context of moment matrix block ranks, however, especially in the pre-asymptotic regime and for typical compression error threshold values, its interpretation should be more nuanced.

In this work, the mutual shadow area is utilized for the study and interpretation of the spectra and ranks of moment matrix blocks, in the asymptotic and pre-asymptotic regimes. 
It is first shown that the mutual shadow area and length (in 3-D and 2-D, respectively), which scale asymptotically as \((ka)^2\) and \(ka\) (as do the leading term in the bound on the number of multipoles), can be used for the accurate prediction of the knee index for interactions with a substantial ``broadside" component.  
The mutual shadow length can also be used for the prediction in 3-D quasi-planar configurations. 
The non-silent remainder of the spectrum, i.e., between the predicted knee index and the crossing of the (greater than machine precision) compression threshold, is shown to be of a width that scales slower with \(ka\), making it asymptotically negligible in determining the rank. 
In some settings, the scaling is as slow as  \(\mathcal{O}(ka)\) and \(\mathcal{O}(1)\), for surfaces in 3- and 2-D, respectively. 

 This behavior is explained through closer examination of the orthonormal bases associated with the two parts of the spectrum and their localization properties in the spatial and Fourier domains.
The findings allude to the association of the spectrum up to the knee with aperture radiation and of the interacting remainder with complementary diffraction phenomena. In some cases, including those common in numerical rank-scaling studies, the remainder can be attributed to end-point contributions or ``edge" diffraction.  
The number of ``aperture DoF" diminishes with the decrease in the aspect angle between the interacting domains and vanishes entirely for quasi-planar configurations. 
This explains the pre-asymptotic skew from the expected \(\mathcal{O}(N)\) rank-scaling observed for moderate-sized geometries and low threshold values.

The remainder of this article is organized as follows: Section~\ref{S:MomentMatrix} provides a brief background on moment matrices and low-rank approximation of their blocks. Section~\ref{S:MutualShadow} defines the mutual shadow area and shadow length for finite source and observer domain, demonstrate their calculation, and shows their relation to moment matrix blocks spectra. Section~\ref{S:Spectrum} analyses and interprets the spectral components with regard to the shadow area. Section~\ref{S:Conclusions} provides additional discussion and directions for follow-up research.

\section{Moment Matrix Block Spectrum and Rank}\label{S:MomentMatrix}

\begin{figure}
    \centering
    \includegraphics[width=\linewidth]{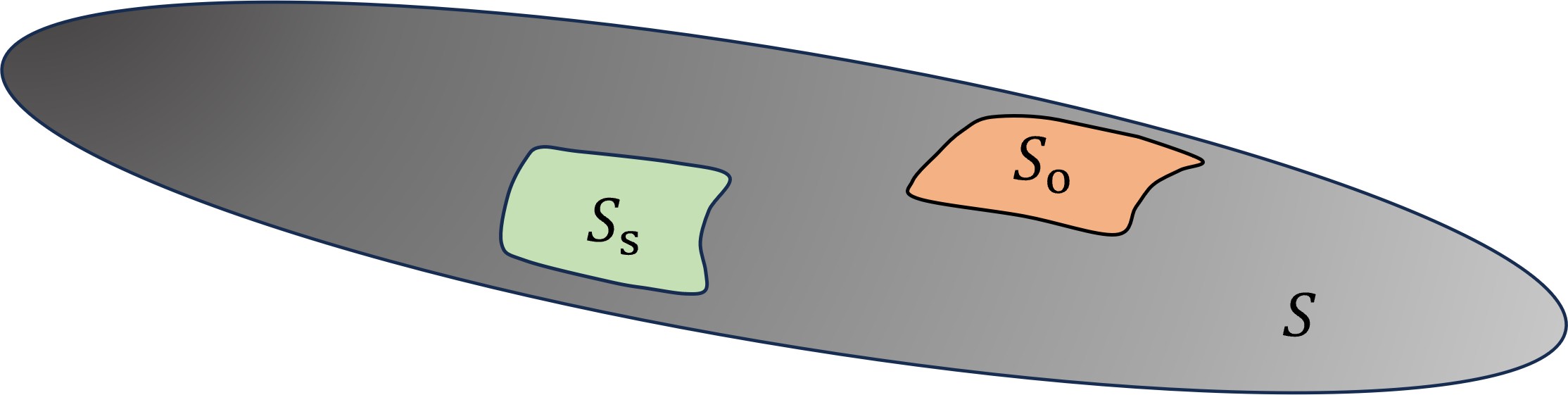}
    \caption{Source $S_\T{s}$ and observer $S_\T{o}$ subdomains of a scatterer \(S\).}
    \label{fig:subdomains_of_S}
\end{figure}

The scattering of waves by an impenetrable scatterer \(S\), see Fig.~\ref{fig:subdomains_of_S}, can be described by using surface IEs. 
For time-harmonic problems (a time dependence $ \eu^{\ju\omega t}$ is assumed and suppressed), the scattered field is presented as an integral on \(S\) (here written, for simplicity, loosely in scalar form) of an unknown surface distribution of sources \(J(\vec{r}'), \vec{r}'\in{S}\) that radiate according to an oscillatory kernel \(K(\vec{r},\vec{r}')\).  Together with the incident field, it satisfies the governing equations outside of \(S\) and the boundary conditions on \(S\).
For a homogeneous background with wavenumber $k$, denoting $R=|\vec{r}-\vec{r}'|$, the scalar kernel is of the form $ K(\vec{r},\vec{r}') =  \mathrm{exp}(-\ju kR)/R $ in 3-D and $ K(\vec{r},\vec{r}') =  H_0^
{(2)}(kR) $ ($H_0^
{(2)}(\bullet)$ being the zeroth order Hankel function of the second kind) in 2-D~\cite{Harrington1968}. 
In its discretization, e.g., by using the MoM with \(N\) basis and testing functions, the surface IE translates to a system of linear equations of the form 
\(\M{Z}\M{i} = \M{v}\)
where \(\M{v}\in \mathbb{C}^N\) and \(\M{i}\in \mathbb{C}^N\) are the tested excitation and unknown coefficient vectors and \(\M{Z}\in \mathbb{C}^{N \times N}\) system or impedance matrix~\cite{Harrington1968}. 

\begin{figure}
    \centering
    \includegraphics[width=\linewidth]{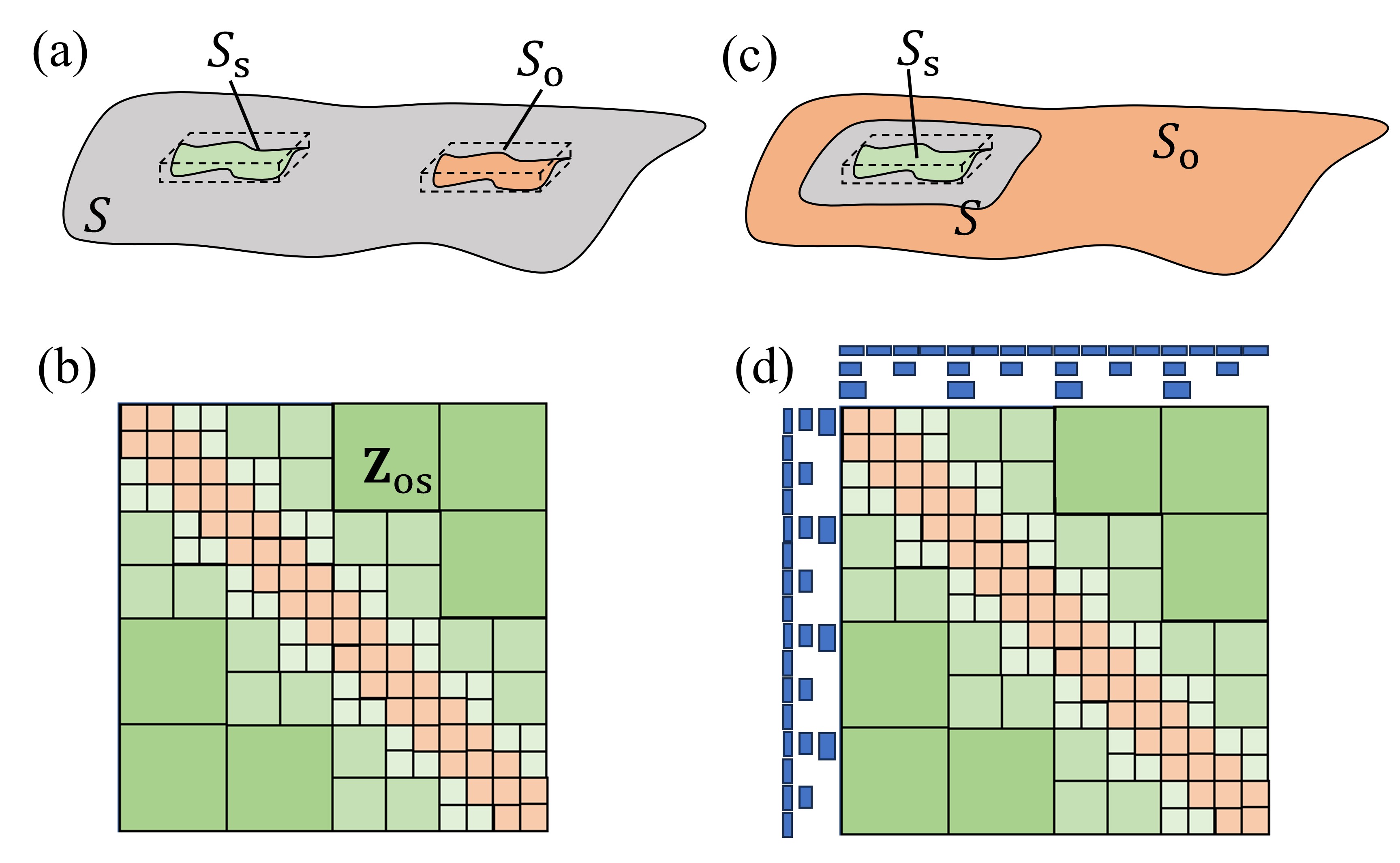}
    \caption{Admissible blocks and corresponding interactions for rank and compressive-bases revealing: (a) Strong admissibility interactions in non-nested bases settings. (b) Non-nested bases hierarchical matrix block structure. (c)  Interaction for strong admissibility nested-bases compression. (d) Nested-bases compression of a hierarchical matrix~\cite{bormDatasparseApproximationNonlocal2007}.}
    \label{fig:H_mat_admissibile}
\end{figure}

Fast solvers that rely on algebraic compression of \(\M{Z}\) use a hierarchical geometric clustering of the basis and testing functions on \(S\) in order to partition \(\M{Z}\) into a hierarchy of blocks. 
In the hierarchy, blocks that are identified as good candidates (or ``admissible") for compression are expressed by a more economical approximate form. 
Blocks that do not satisfy the admissibility criterion are further partitioned. 
This continues until sufficiently small clusters and block sizes are reached. 
The criterion for admissibility is usually based on the sizes of the  subdomains associated with the clusters  and the distance between them. 
Let  \(S_{\mathrm{s}}\) and \(S_{\mathrm{o}}\) be two subdomains resulting from the hierarchical partitioning of \(S\) (as illustrated in Fig. \ref{fig:subdomains_of_S}), on which \(N_{\mathrm{s}}\) and \(N_{\mathrm{o}}\) basis and testing functions are defined. 
These subdomains can be circumscribed by spheres of radii \(R_{\mathrm{s}}\) and \(R_{\mathrm{o}}\) respectively. 
The block that represents the interaction between these clusters is denoted \(\M{Z}_{\mathrm{os}} \in \mathbb{C}^{N_{\mathrm{o}} \times N_{\mathrm{s}}}\). 
The admissibility criterion often takes the form (see\cite{bormDatasparseApproximationNonlocal2007})
\begin{equation}
    R_{\mathrm{s}}+R_{\mathrm{o}}<\eta d,
\end{equation}
where \(d\) denotes the distance (according to some definition) between the subdomains and \(\eta\) is a parameter that determines the separation required for admissibility. 
Strong admissibility criteria identify as compressible only interactions between strictly non-overlapping subdomains.
Let \(\M{Z}_{\mathrm{os}}\) be an admissible block, representing an interaction of the form
\begin{equation}
    \int_{S_\mathrm{s}} K(\vec{r}, \vec{r}') J(\vec{r}') \, \diff S' ,\quad    \vec{r} \in S_\mathrm{o}
    \label{eq:integral_surface}
\end{equation}
In fast solvers, compressed blocks \(\M{Z}_{\mathrm{os}}\) typically correspond to  \(S_{\mathrm{s}}\) and \(S_{\mathrm{o}}\) of the same partitioning level and that are of similar size (see Fig.~\ref{fig:H_mat_admissibile}(a)).
They can be expressed either directly (in a non-nested manner, as is the case in \(\mathcal{H}\) - matrix type solvers , ~\cite{hackbuschSparseMatrixArithmetic1999, zhaoAdaptiveCrossApproximation2005}, see Fig.~\ref{fig:H_mat_admissibile}(b)) or implicitly through multiplication of previously computed compressive nested bases from lower tree levels (as in, for example, the \(\mathcal{H}^2\) - matrix approach ~\cite{bormDatasparseApproximationNonlocal2007}).

When compressed independently, these are the spectra of the blocks \( {\M{Z}}_{\mathrm{os}}\) that directly dictate the computational costs. 
Let \( \M{U}_\mathrm{os} \M{\Sigma}_\mathrm{os}  {\M{V}^{\dagger}_\mathrm{os}} =  \M{Z}_{\mathrm{os}}\) be the full singular value decomposition (SVD) of  \(\M{Z}_{\mathrm{os}}\), where \(\M{U}_\mathrm{os} \in \mathbb{C}^{N_\mathrm{o} \times N_\mathrm{o} }\), \(\M{V}_\mathrm{os} \in \mathbb{C}^{N_\mathrm{s} \times N_\mathrm{s} }\) are unitary matrices and \(\M{\Sigma}_\mathrm{os} \) is a diagonal matrix containing the SVs \(\sigma_n\) in a decreasing order, and the dagger superscript 
denotes Hermitian transpose. 
A rank \(\mathcal{R}^*(\tau)\) approximation \(\M{Z}_{\mathrm{os}} \approx \M{A}\M{B}^\dagger\) (\(\M{A} \in \mathbb{C}^{N_\mathrm{o} \times \mathcal{R}^* }\) \(\M{B} \in \mathbb{C}^{N_\mathrm{s} \times \mathcal{R}^* }\)) is the smallest rank LR approximation of relative error (in the spectral norm) \(\tau\) if \(\mathcal{R}^*\) is the number of singular values satisfying \(\sigma_n/\sigma_1 >\tau\).  
In nested bases compression (see Fig.~\ref{fig:H_mat_admissibile}(d)), the ranks that determine the solver costs are those of interactions between subdomains in each level and entire complementary regions of (typically non-touching) observers. 
Denoting \(\hat{S}_{\mathrm{s}}\) and \(\hat{S}_{\mathrm{o}}\) the non-neighboring parts of the complementary regions to \(S_{\mathrm{s}}\) and \(S_{\mathrm{o}}\) (see Fig. \ref{fig:H_mat_admissibile}(c)), the  ranks of the corresponding compressive bases are those of the blocks denoted \(\M{Z}_{\mathrm{\hat{s} s}}=\M{U} \M{\Sigma} {\M{V}^{\dagger}_\mathrm{s}}\) and \(\M{Z}_{\mathrm{o\hat{o}}} =\M{U}_\mathrm{o} \M{\Sigma} {\M{V}^\dagger}\).

This work focuses on the study of the spectra and ranks of blocks \(\M{Z}_{\mathrm{os}}\), \(\M{Z}_{\mathrm{\hat{s}s}}\), and \(\M{Z}_{\mathrm{o\hat{o}}}\) in strong admissibility interaction solvers, their ranks, and the physical characteristics of the compressive bases that are the first \(\mathcal{R}^*(\tau)\) columns of \(\M{U}_{\mathrm{os}}\),  \(\M{V}_{\mathrm{os}}\),  \(\M{U}_{\mathrm{o}}\) and  \(\M{V}_{\mathrm{s}}\). 
The study relies on the identification of the role of the mutual shadow area and length in the prediction of the asymptotic rank of the interactions, as discussed in the next section.

\section{Mutual Shadow and the Asymptotic Scaling of Moment Matrix Blocks Ranks}\label{S:MutualShadow}
\begin{figure}
    \centering
    \includegraphics[width=\linewidth]{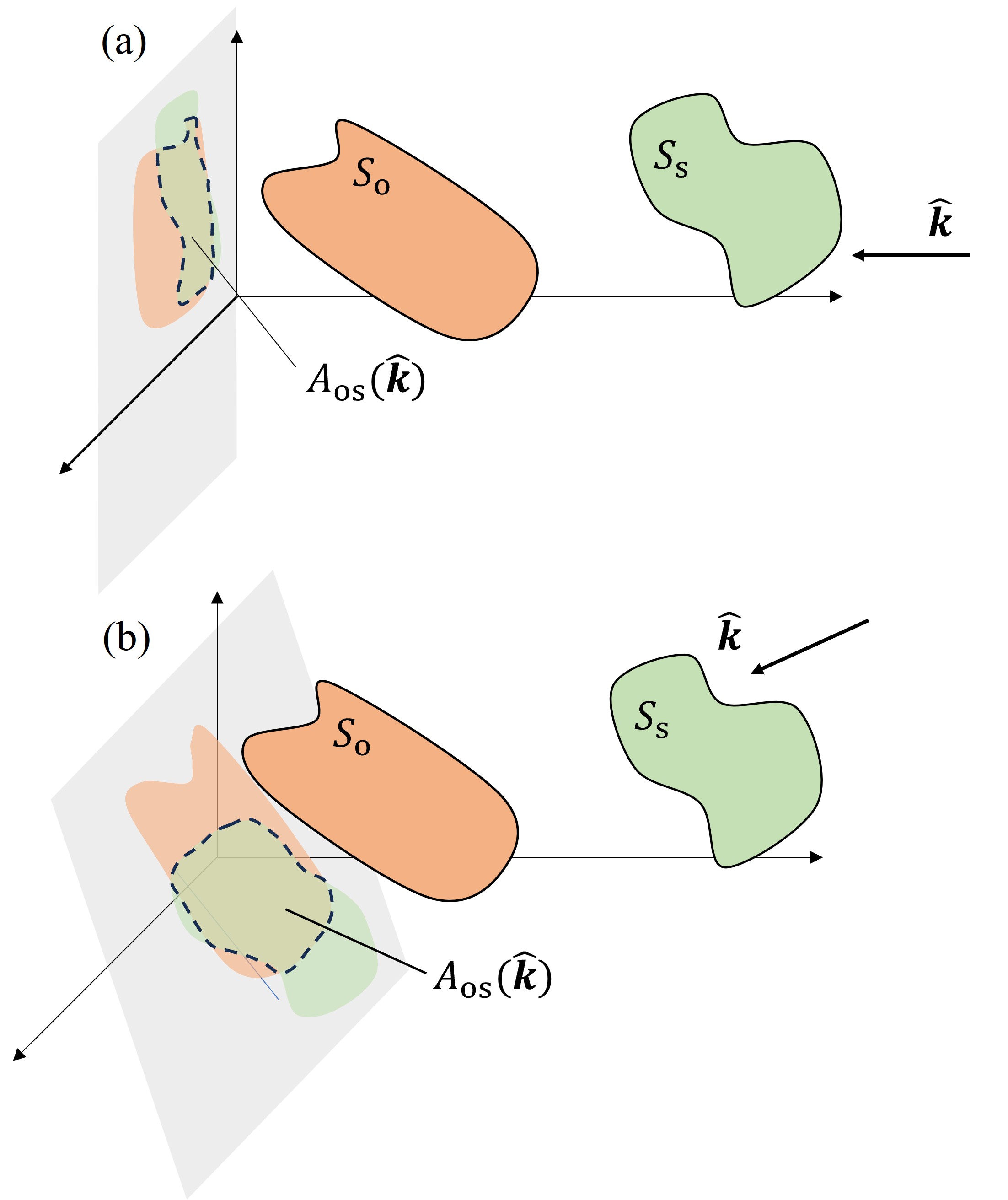}
    \caption{Mutual shadow areas $A_\T{os}(\UV{k})$ for two cases of the illumination direction \(\UV{k}\).}
    \label{fig:mutual_shadow_area_illustrated}
\end{figure}

The total mutual shadow area (in 3-D) and length (in 2-D) were introduced in~\cite{Gustafsson2025c} as geometric quantities that enable accurate prediction of the number of DoF in a MIMO communication channel between arrays of transmitters and receivers. 
Consider a pair of well-separated source and observe regions  \(S_{\mathrm{s}}\) and \(S_{\mathrm{o}}\) to which the receivers and transmitters are confined. 
In our case, these are basis and testing function clusters.  
For a plane wave with a wave-vector \(\UV{k}k\),  approaching \(S_{\mathrm{o}}\) from the general direction of \(S_{\mathrm{s}}\), the mutual shadow area and length, denoted \(A_\T{os}(\UV{k})\) and \(L_\T{os}(\UV{k})\), are defined as the overlapping area and length of the shadows cast by \(S_{\mathrm{s}}\) and \(S_{\mathrm{o}}\) on the plane orthogonal to \(\UV{k}\) (see illustration for two different \(\UV{k}\) values in Figs.
\ref{fig:mutual_shadow_area_illustrated}~(a) and (b)).
It can be thought of as the area of a \(\UV{k}\)-directed tube of rays that cross both \(S_{\mathrm{s}}\) and \(S_{\mathrm{o}}\).
The {\it{total}}  mutual shadow area and length are the integrals of contributions \(\T{d}A_\T{os}(\UV{k})=A_\T{os}(\UV{k})\T{d}\Omega_{\UV{k}}\) and \(\T{d}L_\T{os} (\UV{k})\T{d}\varphi_{\UV{k}}\) by all \(\UV{k}\) pointing generally from  \(S_{\mathrm{o}}\) to \(S_{\mathrm{s}}\), defined as 
\begin{equation}
    \mathcal{A}_\mathrm{os}=
    \int {\T{d}A_\mathrm{os}(\UV{k})=\int {A_\T{os}(\UV{k})\,\T{d}\Omega_{\UV{k}}}}
    \label{eq:area_general}
\end{equation}
and 
\begin{equation}
    \mathcal{L}_\mathrm{os}=
    \int {\T{d}L_\mathrm{os}(\UV{k})=\int {L_\T{os}(\UV{k})\,\T{d}\varphi_{\UV{k}}},}
    \label{eq:length_general}
\end{equation}
respectively. 
For scalar (single-layer potential) kernel interactions, the channel asymptotic number of DoF is approximated in~\cite{Gustafsson2025c} by 
\begin{equation}
    \mathcal{N}_\mathrm{os}^\mathrm{} =\mathcal{A}_\T{os}/\lambda^2
    \label{eq:NAdef}
\end{equation}
and 
\begin{equation}
    \Tilde{\mathcal{N}}_\mathrm{os}^\mathrm{} =\mathcal{L}_\T{os}/\lambda,
    \label{eq:NLdef}
\end{equation}
in 3-D and 2-D, respectively. For vector sources and fields and dyadic IE kernel these numbers should be doubled.  

The general expressions in~\eqref{eq:area_general} and~\eqref{eq:length_general} can be simplified for certain classes of \(S_{\mathrm{o}}\) and \(S_{\mathrm{s}}\) pairs. 
For example, for surface domains \(S_{\mathrm{o}}\) and \(S_{\mathrm{s}}\) where a line-of-sight exists between every pair \(\V{r} \in S_{\mathrm{o}}\) and \(\V{r}' \in S_{\mathrm{o}}\), denoting \(\UV{n}\) and \(\UV{n}'\) the unit normal vectors to these surfaces at \(\V{r}\) and \(\V{r}'\), respectively, \(\mathcal{A}_\mathrm{os}\) can be written, following~\cite{Gustafsson2025c}, as 
\begin{equation}
    \mathcal{A}_\T{os}=
    \int_{S_{\mathrm{s}}}\int_{S_{\mathrm{o}}} \frac{|\UV{n}'\cdot (\V{r}-\V{r}')|\ |\UV{n}\cdot(\V{r}-\V{r}')|}{|\V{r}-\V{r}'|^4}\diff S' \diff S.
    \label{eq:Area_LoS_Surface}
\end{equation}
In 2-D settings, the simplification takes the form of
\begin{equation}
    \mathcal{L}_\T{os}=
    \int_{S_{\mathrm{s}}}\int_{S_{\mathrm{o}}} \frac{|\UV{n}'\cdot (\V{r}-\V{r}')|\ |\UV{n}\cdot(\V{r}-\V{r}')|}{|\V{r}-\V{r}'|^3}\diff l' \diff l.
    \label{eq:Length_LoS_Surface}
\end{equation}
These simplified integral expressions, related also to radiation view factors in the context of radiative heat transfer~\cite{Welty+etal2020,Howell+Menguc2011}, can be evaluated numerically, through various choices of surface discretizations and quadrature rules. 
Notably, they are independent of \(k\) such that the cost of numerical integration is of \({\mathcal O}(1)\).
In this work, numerical evaluation of \eqref{eq:Area_LoS_Surface} and \eqref{eq:Length_LoS_Surface} is used almost exclusively. 
In the special rotation-symmetric case of parallel circular disc domains of radius \(a\) and height separation \(d\) (as in the inset of Fig. 4), \eqref{eq:Area_LoS_Surface} can be calculated in closed form as~\cite{Gustafsson2025c} 
\begin{equation}    
	\mathcal{A}_\T{os} = \frac{\pi^2}{4} \big(\sqrt{4a^2+d^2}-d\big)^2.
    \label{eq:Area_discs}
\end{equation}

 \begin{figure}
    \centering
    \includegraphics[width=\linewidth]{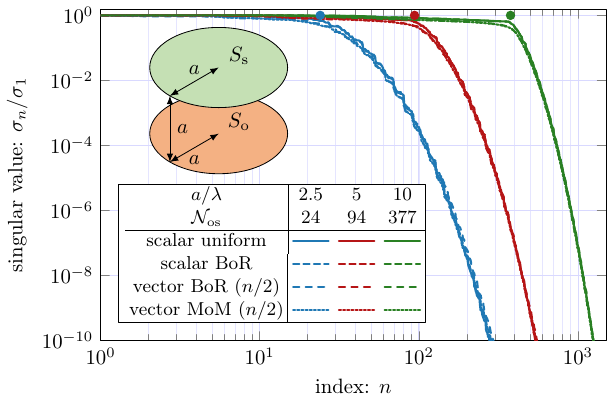}
    \caption{Normalized SVs for two discs evaluated using four different models and three electrical sizes $a/\lambda\in\{2.5,5,10\}$. The asymptotic DoF~\eqref{eq:NAdef} based on~\eqref{eq:Area_discs} is given in the legend.}
    \label{fig:discs_methods}
\end{figure}

The calculation is demonstrated in the following examples, for various geometric settings. 
Consider first the interaction between the two circular discs depicted in the inset of Fig.~\ref{fig:discs_methods}. 
of radii $a \in\{2.5,5,10\}\lambda$ separated by a distance $d = a$, the SVs and $\mathcal{N}_\mathrm{os}$  were computed for the following cases and methods: 
1. single-layer potential interaction between uniformly spaced point sources and observers, 2. single-layer potential interaction between circularly [body of revolution (BoR)] sampled point sources and observers, 3. surface current to tangent electric-field interaction between circularly sampled point sources and observers, and 4. surface current to tangent electric-field interaction with method of moments discretization (rooftop basis and testing functions on uniform polar meshes). 
In Fig.~\ref{fig:discs_methods}, the normalized singular values  \(\sigma_n/\sigma_1\) are presented against the index \(n\) for the single layer interactions. 
For the vector interactions, \(\sigma_{n/2}/\sigma_1\) are presented against \(n\).
For the uniformly and circularly sampled currents and fields, a randomized SVD~\cite{Halko+etal2011, martinssonRandomizedNumericalLinear2020} computation was used, with fast Fourier transform (FFT)-based acceleration of the matrix-vector multiplications.
 A very good agreement is observed between the SV curves computed for all four cases.
 This reaffirms the common notion that a point-based model is representative of the spectral analysis of strong-admissibility interactions between sets of basis and testing functions.  
 Thus, moving forward, single-layer potential interactions between point-source and -observer sets are considered. 
 In the SV curves, the expression in \eqref{eq:Area_discs} provides an accurate prediction of the knee. 
 It is already worth noting also that, depending on \(\tau\), \(\mathcal{R} (\tau)\) can be significantly greater (in this case, by up to half an order of magnitude) than \(\mathcal{N}_\mathrm{os}^\mathrm{}\).
 
\begin{figure}
    \centering

    \includegraphics[width=\linewidth]{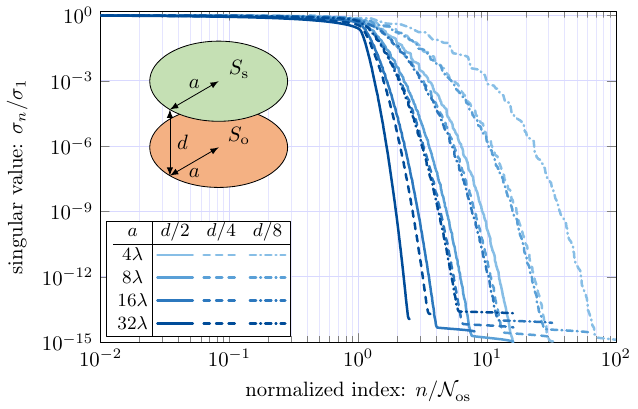}    
    \caption{Normalized SVs on scaled axis for parallel circular discs of various configurations of \(d/a\)  and \(a/\lambda\). }
    \label{fig:discs_SVs_scaling}
\end{figure}

Next, the scaling of the knee-point location with size is examined using the example in the inset of Fig.~\ref{fig:discs_SVs_scaling}. 
This is done for three geometric settings, \(d \in \{2a,4a, 8a\} \), and ranges of \(a\). 
For these configurations, \(\mathcal{N}_\mathrm{os}^\mathrm{}\) (computed using \eqref{eq:Area_discs}) ranges from $\mathcal{N}_\mathrm{os}^\mathrm{}\approx 2.4$ for \(d = 8a = 32\lambda\) to  $\mathcal{N}_\mathrm{os}^\mathrm{} \approx 1734$ for  \(d = 2a = 64\lambda\).
Here, \(\sigma_n/\sigma_1\) are plotted versus a normalized index \(n/\mathcal{N}_\mathrm{os}^\mathrm{}\). 
When the knee of the SVs curve shows at \(n/\mathcal{N}_\mathrm{os}^\mathrm{}=1\), the prediction is said to be accurate. 
In Fig.~\ref{fig:discs_SVs_scaling}, the knee is more identifiable and predicted with increasing accuracy as \(a/\lambda\) grows. 
In the remainder of the spectrum, beyond the predicted knee index, the slope becomes steeper with  \(a/\lambda\), making \(\mathcal{N}_\mathrm{os}^\mathrm{}\) an increasingly more accurate estimate of the rank, even for lower \(\tau\) values. 
Comparing the different geometric configurations, the greater \(d/a\) the later with \(a/\lambda\) the convergence to a step function. 
In other words, \(\mathcal{N}_\mathrm{os}^\mathrm{}\) predicts the rank more accurately for greater  \(a/\lambda\) and \(a/d\) values.

\begin{figure}
    \centering
    
    \includegraphics[width=\linewidth]{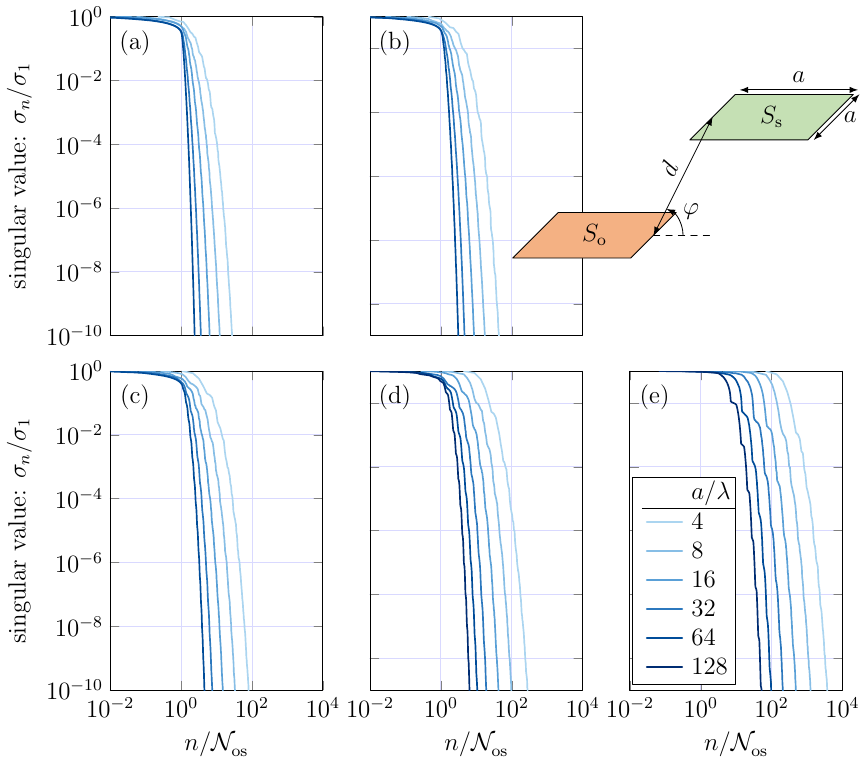}        
    \caption{Normalized SVs on scaled axis for slanted parallel plates for various \(a/\lambda\). (a) \(\varphi\) = \(\pi/2\), (b) \(\varphi\) = \(3\pi/8\), (c) \(\varphi\) = \(\pi/4\), (d) \(\varphi\) = \(\pi/8\), (e) \(\varphi\) = \(0.1\).}
    \label{fig:squares_slanted_SVs_scaling}
\end{figure}

The computation of SV curves is repeated (this time using the fast spectrum revealing method in ~\cite{kelleyIterativeRandomSampling2023}) for squared plate domains with edge length \(a\), shifted both horizontally and vertically with respect to each other, as in the Fig. \ref{fig:squares_slanted_SVs_scaling} inset. 
Here, the distance between the closest edges of the two squares is \(d= a\) and the angle \(\varphi\) determines the horizontal and vertical shifts between the plates. 
The knee index prediction \(\mathcal{N}_\mathrm{os}^\mathrm{}\) is accurate for large \(a/ \lambda\) for a broad range of \(\varphi\) values. 
However, for small \(\varphi\) values, slower convergence to a step-like function is observed.  
For the smallest \(\varphi\) value examined, as the geometric configuration approaches quasi-planar, the knee was inaccurately predicted in the entire examined range of \(a/\lambda\) values. 
These results are in agreement with those in  Fig.~\ref{fig:discs_SVs_scaling}: for a given \(ka\), when the aspect angle between the source and observer domains becomes sufficiently small, due to either the distance or the offset, \(\mathcal{N}_\mathrm{os}^\mathrm{}\) reduces until it no longer predicts the knee correctly.

\begin{figure}
    \centering
    \includegraphics[width=0.95\linewidth]{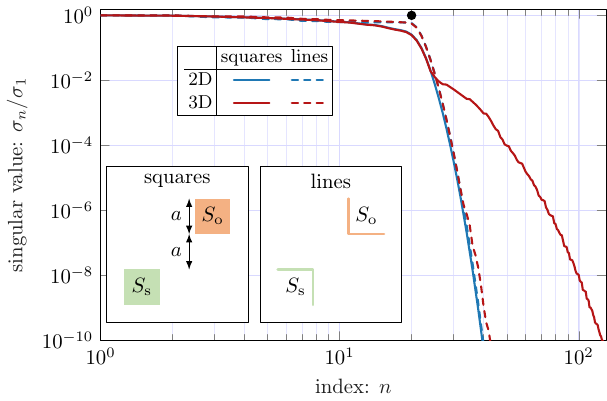}
    \caption{Normalized SVs for planar geometries in 2-D (blue) and 3-D (red). Squares are shown with solid curves and line geometries with dashed curves, for $a=30\lambda$. Marker at $n=20\approx\Tilde{\mathcal{N}}$ based on mutual shadow length of two lines $\mathcal{L}_\T{os}\approx 0.67a$ using~\eqref{eq:Length_LoS_Surface}.}
    \label{fig:SVs_2D3D_comp}
\end{figure}

In the extreme case of (quasi-) planar geometries, \(\mathcal{N}_\mathrm{os}^\mathrm{}\) may vanish entirely. 
In the next example, we examine whether, for such scenarios, the mutual shadow {\it length}-based \(\Tilde{\mathcal{N}}_\mathrm{os}^\mathrm{}\) in~\eqref{eq:NLdef} is a more suitable predictor of the knee index. 
To this end, the SVs are computed for squared domains of edge length \(a=30\lambda\), placed in the plane and shifted diagonally, by \(2a)\) in both the \(x\) and \(y\) directions. 
This is done for both the 3-D and 2-D scalar kernels.
To the results, in Fig~\ref{fig:SVs_2D3D_comp}, a marker that indicates \(\Tilde{\mathcal{N}}_\mathrm{os}^\mathrm{}\), computed using \eqref{eq:Length_LoS_Surface} in-2D, is added. 
In solid lines, it can be seen that while the SV curve for the 3-D kernel exhibits a more complicated behavior after the knee, the 3-D and 2-D kernel SV curves match well until a little beyond the knee, which is predicted accurately by  \(\Tilde{\mathcal{N}}_\mathrm{os}^\mathrm{}=20\).
For comparison, dashed SV curves are plotted for a geometry of the same cumulative mutual shadow, composed of the front corner edges of the squared domains.
These SV curves exhibit the simpler, step-like, behavior and agree well until very low values.
Potential explanations for these differences will be discussed later.
Here, it is deduced that the mutual shadow length prediction of the knee is relevant for planar settings.

\begin{figure}
    \centering
    \includegraphics[width=0.95\linewidth]{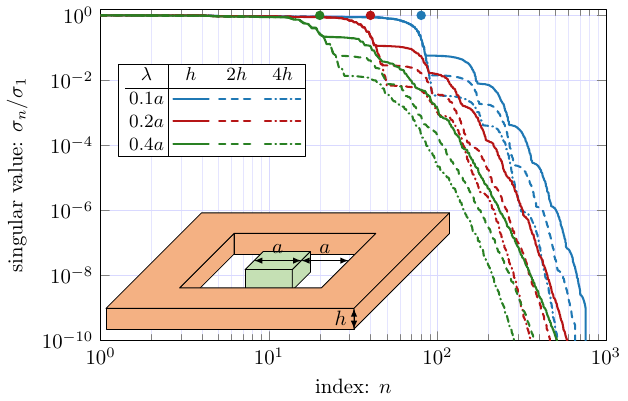}
    \caption{\(\sigma_n /\sigma_1\) for interactions typical for nested-bases compression settings in a quasi-planar geometric configuration.}
    \label{fig:SVs_Slab}
\end{figure}

Lastly, a {\it quasi}-planar case is examined. 
This is done for strong admissibility interactions that are typically encountered in nested bases fast solver settings and represented by \(\M{Z}_{\mathrm{\hat{s} s}}=\M{U} \M{\Sigma} {\M{V}^{\dagger}_\mathrm{s}}\).
Consider the configuration in the inset of Fig.~\ref{fig:SVs_Slab}. 
The normalized SVs are plotted, for various values of the edge length  \(a/\lambda\) and a subwavelength structure height \(h\in\{1,1/2,1/4\}\lambda\) against the index \(n\). Also indicated are the predicted knee indices \(\Tilde{\mathcal{N}}_\mathrm{os}^\mathrm{}\), computed using \eqref{eq:Length_LoS_Surface} for the 2-D cross-section of the geometry. 
Here, the knee indices are predicted well, already for relatively small \(a/\lambda\) values and with increasing accuracy as \(a/\lambda\) grows. 
For each $a/\lambda$ case, they are identical for all $h/\lambda$ values.
Beyond \(\Tilde{\mathcal{N}}_\mathrm{os}^\mathrm{}\), the spectrum shape depends on $h/\lambda$.

To conclude: The mutual shadow area, with  \((ka)^2\) scaling, provides a good prediction  \(\mathcal{N}_\mathrm{os}^\mathrm{}\) of the knee index for 3-D interactions with a sufficiently large aspect angle between the source and observer subdomains. 
For these geometries, it is also an asymptotically accurate predictor of the rank.
For pre-asymptotic sizes, a correction for the prediction of the knee is needed.
For quasi-planar configurations, where the shadow area vanishes, the correction, in the form of \(\Tilde{\mathcal{N}}_\mathrm{os}^\mathrm{}\), based on the mutual shadow {\it length} for the 2-D  cross-section geometry, is in its entirety a good predictor of the knee index. 
There, it is also an asymptotically accurate predictor of the rank and scales as \(ka\).

\section{Mutual Shadow and the Moment Matrix Block Rank and Spectrum}\label{S:Spectrum}

  So far, it was shown that \(\mathcal{N}_\mathrm{os}^\mathrm{}\) is an asymptotically accurate estimate of \(\mathcal{R}\). 
  In that, it resembles the leading  \(\propto(ka)^2\) term in the bound on the number of spherical functions (or multipoles) required for describing the fields produced by a finite domain ~\cite{Colton+Kress1992}. 
  The corrections to this term scale slower with \(ka\) and depend on the prescribed  accuracy for the representation. Here too, for practical (and reasonably large) domain sizes, \(\mathcal{R} (\tau) > \mathcal{N}_\mathrm{os}^\mathrm{} \) depends greatly on  \(\tau\) in a manner that is determined by the geometric configuration and size. 
  In what follows, the spectrum and rank are examined with respect to \(\mathcal{N}_\mathrm{os}^\mathrm{}\) in that pre-asymptotic regime.

  \begin{figure}
    \centering
    \includegraphics[width=0.95\linewidth]{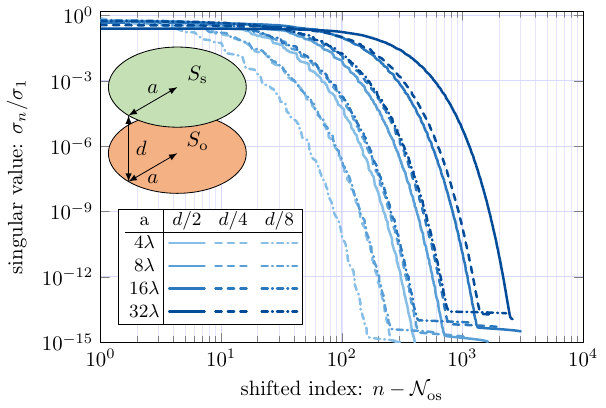}    
    \caption{Remainder parts of SV curves for interactions between the discs in Fig. \ref{fig:discs_SVs_scaling}.}
    \label{fig:discs_SVs_scaling_tails}
\end{figure}

\begin{figure}
    \centering
    \includegraphics[width=0.95\linewidth]{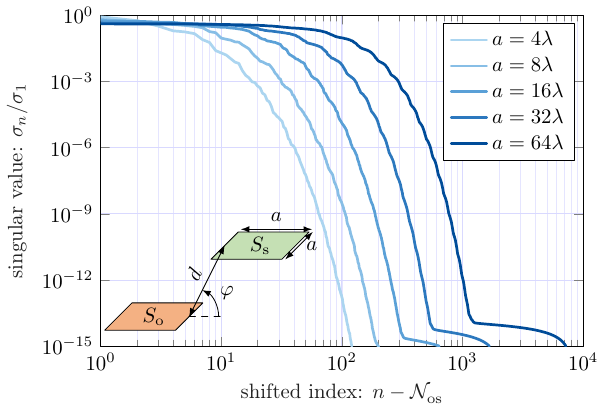}    
    \caption{Remainder parts of SV curves for interactions between the squares in Fig.~\ref{fig:squares_slanted_SVs_scaling} (c).}
    \label{fig:squares_slanted_SVs_scaling_tails}
\end{figure}

 Consider the examples in  Fig.~\ref{fig:discs_SVs_scaling} and let us focus on the non-silent (or interacting) remainder of the spectrum, i.e., on the indices greater than \(\mathcal{N}_\mathrm{os}^\mathrm{}\).  
Fig.~\ref{fig:discs_SVs_scaling_tails} plots \(\sigma_n /\sigma_1\) against a shifted (but not scaled) axis \((n-\mathcal{N}_\mathrm{os}^\mathrm{})\).
The remainder curves exhibit a knee that shifts with \(a/\lambda\), beyond which the curves drop in parallel. 
 That is, the correction to the prediction \(\mathcal{N}_\mathrm{os}^\mathrm{}\) for the knee index scales roughly linearly with \(ka\).
 It is unclear, however, whether it can also be calculated using simple geometric considerations.
 Asymptotically, \(\mathcal{R} (\tau) - \mathcal{N}_\mathrm{os}^\mathrm{} \) is also, roughly, doubled with the doubling of \(a/\lambda\). 
 This contribution to \(\mathcal{R}\), which is, naturally larger for lower \(\tau\) values, is more significant for small \(a/\lambda\) values.  
 Growing slower than \( \mathcal{N}_\mathrm{os}^\mathrm{} \) with \(a/\lambda\), it is {\it asymptotically} negligible.
  For the geometry in Fig.~\ref{fig:squares_slanted_SVs_scaling} and the case in Fig. \ref{fig:squares_slanted_SVs_scaling} (c) (\(\varphi = \pi/4\)), Fig. \ref{fig:squares_slanted_SVs_scaling_tails} shows the SVs for indices greater than \(\mathcal{N}_\mathrm{os}^\mathrm{}\).
    A similar behavior is observed, with scaling that, although slightly faster than \(ka\), is significantly slower then \((ka)^2\).

\begin{figure}
    \centering
    \begin{tikzpicture} 
    \node at (0,0) {\includegraphics[width=0.95\linewidth]{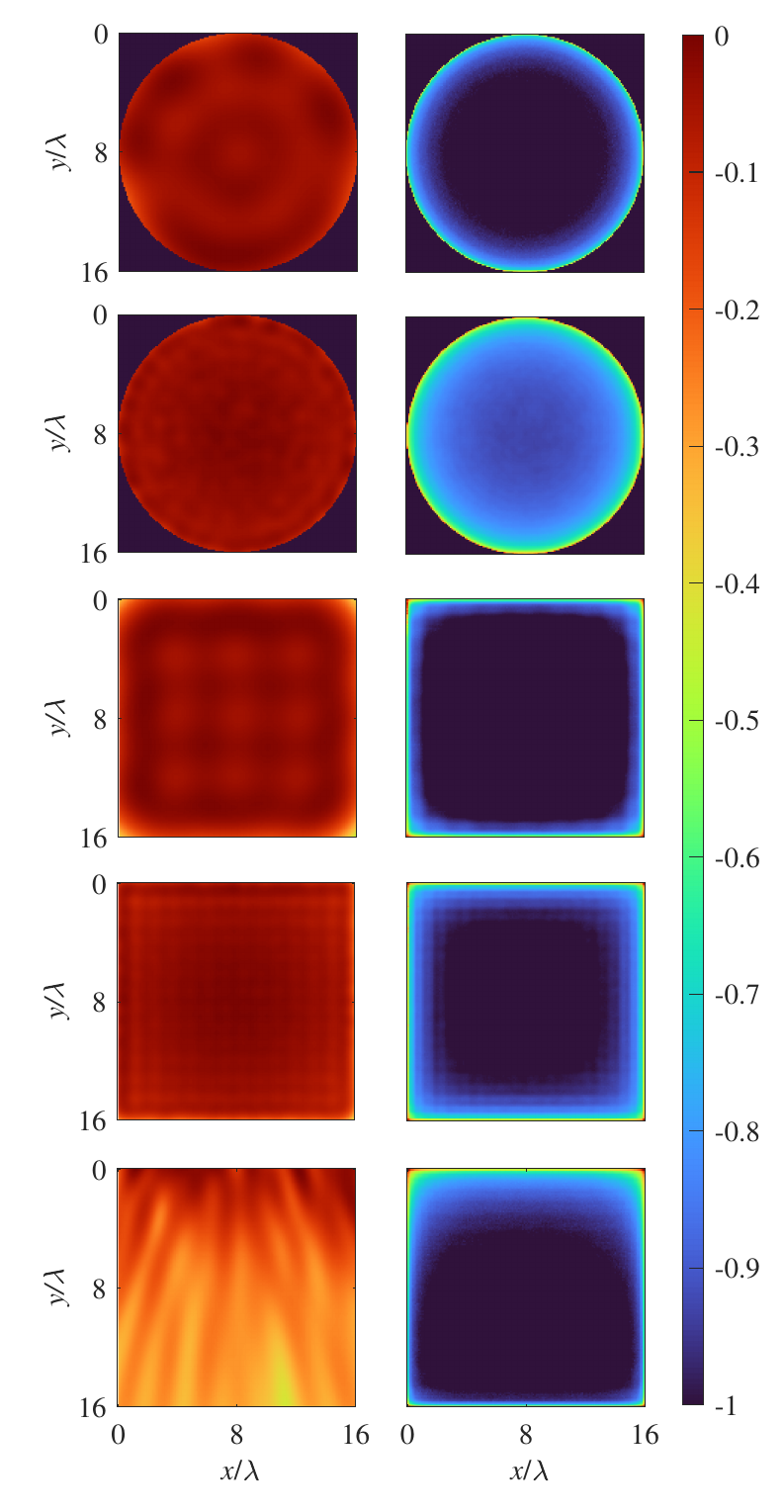}};
    \node at (-4,8) {(a)};
    \node at (-4,4.9) {(b)};
    \node at (-4,1.8) {(c)};
    \node at (-4,-1.3) {(d)};
    \node at (-4,-4.4) {(e)};
    \end{tikzpicture}
    \caption{Logarithm of local average of \( \M{U}_\mathrm{os} \) columns absolute values for (left) first \( \mathcal{N}_\mathrm{os}^\mathrm{} \) and (right) the slope \(\mathcal{R} (\tau) - \mathcal{N}_\mathrm{os}^\mathrm{} \) indexed singular values. (a) Circular disks \(a/ \lambda = 8\), \(d/a = 8\). (b) Circular disks \(a/ \lambda = 8\). (c) Parallel squared plates \(d/a = 4\), \(a/\lambda = 16\). (d) Parallel squared plates \(d/a = 1\), \(a/\lambda = 16\). (e) Slanted parallel squared plates as in Fig.~\ref{fig:squares_slanted_SVs_scaling}  with \(a/\lambda = 16\), \( \varphi = \pi/4\).}
    \label{fig:modes_plates_avg_abs}
\end{figure}
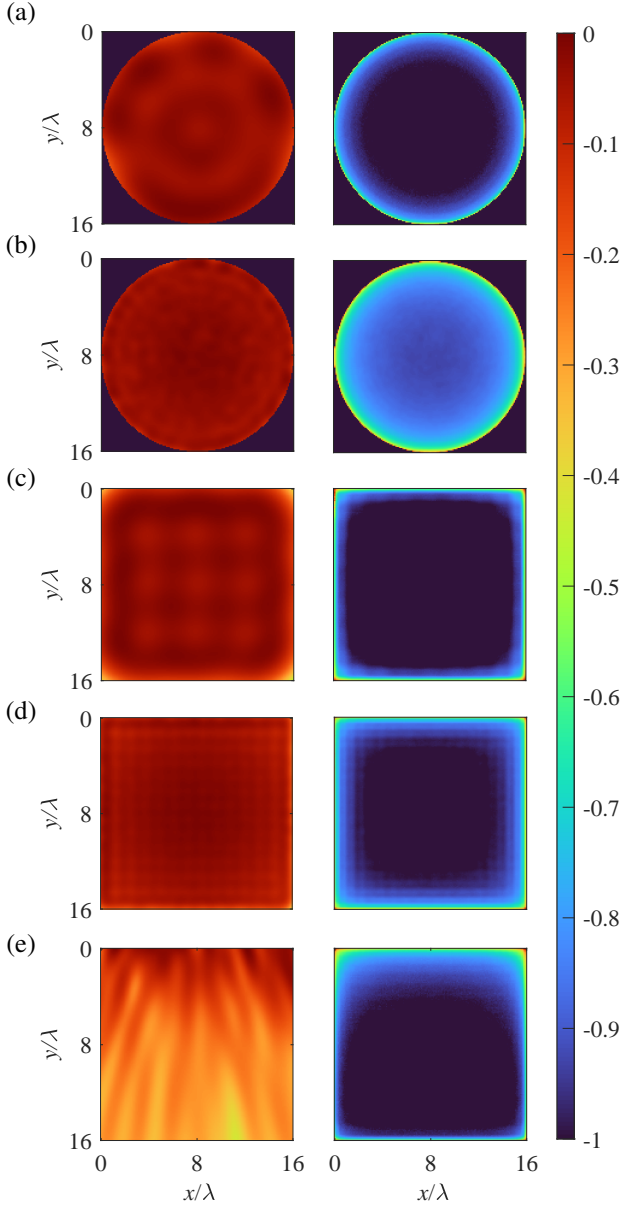

The (close to) linear \(ka\)-scaling in Fig.~\ref{fig:squares_slanted_SVs_scaling_tails} of \(\mathcal{R} (\tau) - \mathcal{N}_\mathrm{os}^\mathrm{} \) may hint that the non-zero (to machine precision) remainder of the SVs curve can be attributed to phenomena associated with the length of the edges of the interacting subdomains, rather than their area. 
This is shown, in Fig. \ref{fig:modes_plates_avg_abs}, for selected examples: pairs of parallel circular disks (Fig.~\ref{fig:discs_SVs_scaling}), with radii \(a= 8 \lambda \) and spacings \(d = 8a\) and \(d = 2a \), squared plates with edge length \(a= 16 \lambda \) and spacings \(d = 4a\)  and \(d = 2a \), and slanted squared plates (Fig.~\ref{fig:squares_slanted_SVs_scaling}) with  \(a = 16\lambda\), \( \varphi = \pi/4\). 
For these configurations, using the randomized range revealing method in \cite{Halko+etal2011,martinssonRandomizedNumericalLinear2020}, the root mean squares along the row dimension of \( \M{U}_\mathrm{os} \) are computed, for the first \( \mathcal{N}_\mathrm{os}^\mathrm{} \) columns and the next \(\mathcal{R} (\tau) - \mathcal{N}_\mathrm{os}^\mathrm{} \) ones (here, for \(\tau = 10^{-15}\)). 
The results in Fig.~\ref{fig:modes_plates_avg_abs}, for each case, show the common logarithm of the average values, normalized to the maximum in each such map, within a dynamic range of one order of magnitude.
In each of the cases, the averages for the first \( \mathcal{N}_\mathrm{os}^\mathrm{} \) columns shows to be rather uniform in magnitude across the subdomain surface, with mild reduction near the edges. 
For the slanted plates case in Fig.~\ref{fig:modes_plates_avg_abs} (e), higher intensity is on the observer domain is seen on the side that is closer to the source domain. 
For the additional \(\mathcal{R} (\tau) - \mathcal{N}_\mathrm{os}^\mathrm{} \) columns, the average exhibits distinct non-uniformity, with higher intensity near the edges. 
This effect is especially pronounced for the cases of greater separation in Figs. \ref{fig:modes_plates_avg_abs} (a) and (c).  
It is deduced that, for interactions between flat finite surfaces, the \(\mathcal{R} (\tau) - \mathcal{N}_\mathrm{os}^\mathrm{} \)-dimensional remainder interacting subspace is spatially more localized near the subdomain fringes. 
For the examined cases, its contribution to the interaction can be associated with edge diffraction that complements the rank-\( \mathcal{N}_\mathrm{os}^\mathrm{} \) ``aperture" subspace of the interaction.

\begin{figure}
    \centering
    \includegraphics[width=\linewidth]{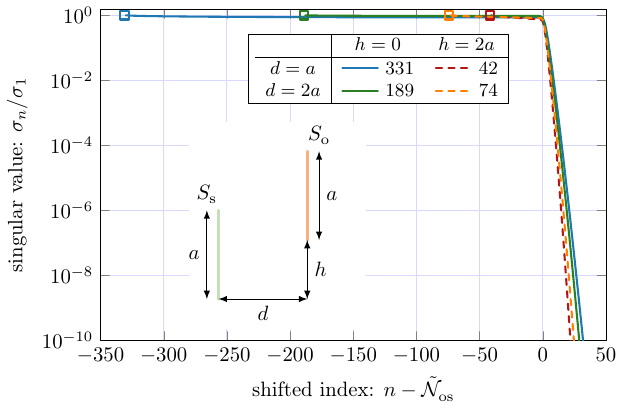}
    \caption{\(\sigma_n /\sigma_1\) for various configurations of the geometry in the inset, varying in \(d/a\) and \(h/a\), according to the legend, with \(a=400
\lambda\). Square markers indicate \(n=1\) for each of the cases. The knee index, as predicted by~\eqref{eq:NLdef} based on~\eqref{eq:ShadowLengthTwoLines}, is reported in the legend.}
    \label{fig:SVs_lines_2D}
\end{figure}

\begin{figure}
    \centering
    \includegraphics[width=0.95\linewidth]{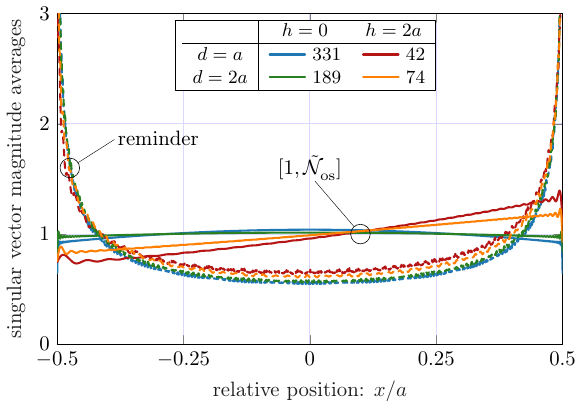}
    \caption{Local mean squares of $\M{U}_\T{os}$ columns absolute values for the cases in Fig.~\ref{fig:SVs_lines_2D}. Solid curves correspond to the first $\Tilde{\mathcal{N}}_\mathrm{os}^\mathrm{}$ columns.  Dashed curves correspond to the next 25 vectors. }
    \label{fig:modes_lines_avg_abs}
\end{figure}

Additional insight into this behavior can be gained from the analysis of 2-D interactions. 
To this end, let us first repeat the experiment of Fig.~\ref{fig:squares_slanted_SVs_scaling_tails} and Fig.~\ref{fig:modes_plates_avg_abs} for configurations of parallel line subdomains of length \(a\) each, spaced and shifted by distances \(d\) and  \(h\), respectively, as shown in the  Fig.~\ref{fig:SVs_lines_2D} inset. 
For this case, the total mutual shadow length in~\eqref{eq:Length_LoS_Surface} reduces to 
\begin{equation}
    \mathcal{L}_\T{os} = \sqrt{d^2+(a+h)^2} - 2\sqrt{d^2+h^2} + \sqrt{d^2+(a-h)^2}
    \label{eq:ShadowLengthTwoLines}
\end{equation}
In Fig.~\ref{fig:SVs_lines_2D}, \(\sigma_n /\sigma_1\) are presented versus \(n - \Tilde{\mathcal{N}}_\mathrm{os}^\mathrm{}\). The accurate prediction of the knee indices is seen at zero. 
The slope following the knee is very steep, making   \(\Tilde{\mathcal{N}}_\mathrm{os}^\mathrm{}\) a good predictor also for the effective rank. 
Here, the slope depends very weakly (or hardly at all) on \(d/a\) and \(h/a\). 
This behavior agrees with that observed for planar apertures in 3-D: the dimension of the interacting remainder subspace beyond \( \Tilde{\mathcal{N}}_\mathrm{os}^\mathrm{}\) scales slower with \(ka\) than \( \Tilde{\mathcal{N}}\) (almost as \(O(1)\)).  
It is next examined whether this subspace can also be attributed to the edges of the lines. 
An exercise similar to that in Fig.~\ref{fig:modes_plates_avg_abs} is carried out, for the geometric configurations of Fig. \ref{fig:SVs_lines_2D}. In Fig.~\ref{fig:modes_lines_avg_abs}, the (not normalized) local mean squares of the absolute values of the singular vectors are plotted. 
 The first \(\Tilde{\mathcal{N}}_\mathrm{os}^\mathrm{}\) basis vectors (solid lines) have, on average, a rather uniform intensity along the domain. The averages on the next 25 basis vectors (dashed lines), which are within the interacting remainder part of the spectrum, show higher intensity near the line edges. 
A tilt, similar to that observed in Fig.~\ref{fig:modes_plates_avg_abs}(e), shows here for the cases of shifted lines.

\begin{figure}
    {\centering
    \includegraphics[width=\linewidth]{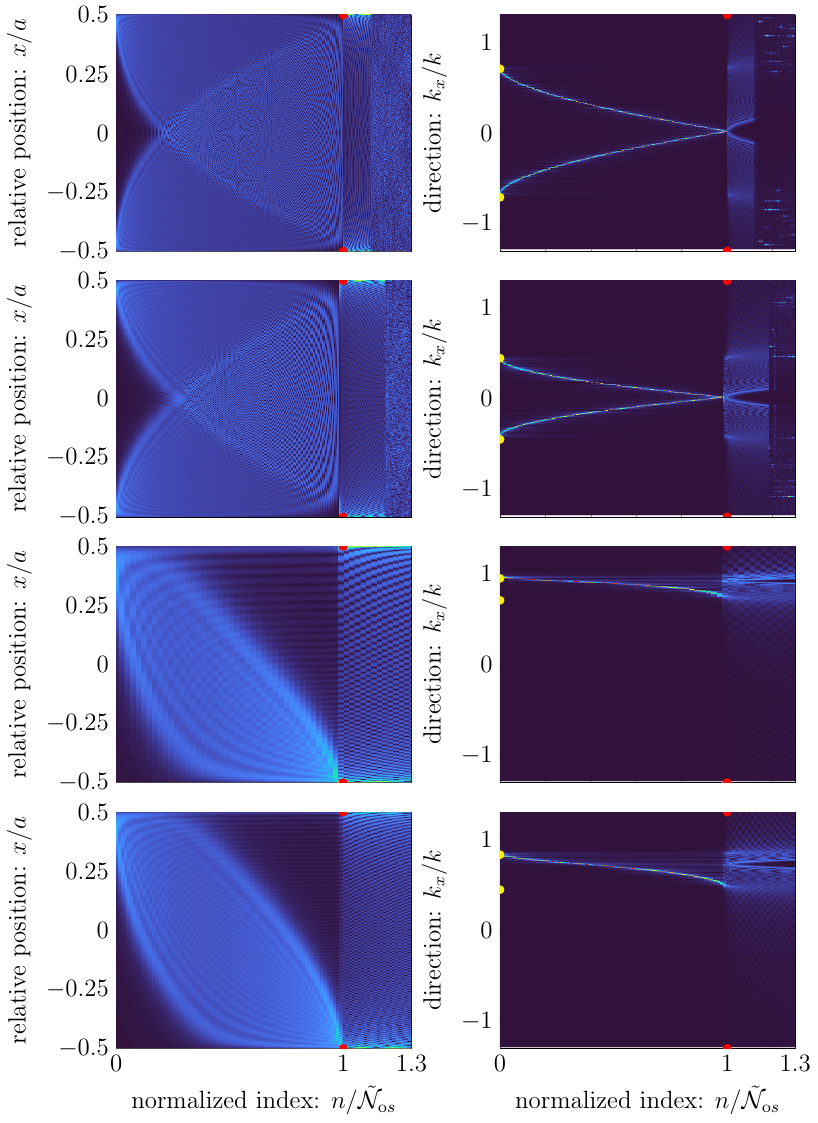}
    \par}
    
    \caption{(left) Singular vectors with their (right) discrete Fourier transforms, for the configurations in Fig.~\ref{fig:SVs_lines_2D}. From top to bottom $(h/a,d /a)\in\{(0,1),(0,2),(2,1),(2,2)\}$. }
    \label{fig:linege}
\end{figure}

\begin{figure}
    \centering
    \includegraphics[width=0.95\linewidth]{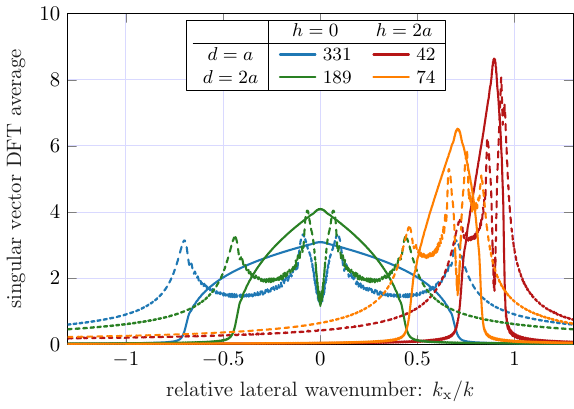}
    \caption{Local mean squares of the discrete Fourier transform of the columns of $\M{U}_\T{os}$ used in Fig.~\ref{fig:SVs_lines_2D}. Solid curves correspond to the first $\Tilde{\mathcal{N}}_\mathrm{os}^\mathrm{}$ columns.  Dashed curves correspond to the next 25 vectors.}
    \label{fig:modes_lines_avg_abs_fft}
\end{figure}

The 2-D case allows us to examine the interacting bases from up-close. 
In  Fig.~\ref{fig:linege}, for the selected geometric configurations in Fig. \ref{fig:modes_lines_avg_abs}, the absolute values of the entries of the first \(1.3~\Tilde{\mathcal{N}}_\mathrm{os}^\mathrm{}\) columns of \( \M{U}_\mathrm{os} \)  are presented. 
In each sub-figure, a red markers near the index axis indicates \(\Tilde{\mathcal{N}}_\mathrm{os}^\mathrm{}\). 
In Figs~\ref{fig:linege}(left), it can be seen that, with surprising precision, \(\Tilde{\mathcal{N}}_\mathrm{os}^\mathrm{}\) marks the partitioning of the domain of the matrix block into two subspaces of distinct wave characteristics. 
The first (approximately  \(\Tilde{\mathcal{N}}_\mathrm{os}^\mathrm{}\)) basis vectors are spatially oscillatory and of a rather uniformly distributed intensity along the line. 
Of the remaining vectors, the ones associated with \(\sigma_n /\sigma_1\) that are greater than machine precision have a concentrated intensity near the edges 
(the remaining vectors are noisy and are in the silent (non-interacting) subspace of the block). 
Interestingly, the discrete Fourier transforms (DFTs) (Figs~\ref{fig:linege}(right)) (i.e., their \(k_x\)-spectral contents) exhibit, for the first vectors, a localized (spectral) behavior, scanning in various directions. 
The scanning is symmetric for \(h=0\) and single-sided for \(h=2a\). The ranges for the scanning lateral wavenumber $k_x$-values correspond to the ranges of directions between the pairs of points in $S_\T{s}$ and $S_\T{o}$ in each configuration and are indicated by yellow markers,  with narrower ranges for the larger \(d/a\) value.  
The next interacting singular vectors are transformed into spectrally non-localized ones, as can be expected for the fields produced by spatially highly-localized sources.
The  magnitude mean squares for the basis vector DFTs, similar to those in Fig.~\ref{fig:modes_lines_avg_abs}, are presented in Fig.~\ref{fig:modes_lines_avg_abs_fft}.

\section{Conclusion and Discussion}\label{S:Conclusions}

The concept of mutual shadow area and length was shown to be useful for interpreting the SV spectrum of moment matrix blocks. 
The examples show that the cumulative mutual shadow area, which scales as \((ka)^2\),
provides a very good prediction of the SV curve knee index for broad ranges of geometric configurations. 
The accuracy of the prediction deteriorates for small aspect angles formed between the source and observer regions.  
The correction to the shadow area-based knee prediction was shown to be slower scaling with \(ka\) and dependent on the geometry.
For cases in which the shadow area vanishes, the correction governs the knee index. 
For quasi-planar cases, it was shown to be predictable using the cumulative mutual shadow length of a 2-D equivalent of the geometry.

The significance of the knee index for moment matrix blocks differs from that 
in communication channels capacity estimation applications. 
Fast solvers require knowing the block rank for securing a prescribed LR approximation relative error level.
Asymptotically, the rank approaches the knee index. 
However, for arbitrary source-observer configuration and sizes, its value depends on the remainder of the spectrum beyond the predicted knee index (which  scale as \((ka)^2\) in 3-D or \(ka\) for quasi-planar and 2-D configurations). 
The remainder of the interacting (to the prescribed threshold) spectrum is of a width that scales much slower with \(ka\).
Furthermore, its scalings is often slower than that of the first correction term to the bound on the required number of multipoles.
For example, for flat interacting domains, it can be as slow as \(O(ka)\) and \(O(1)\) in 3- and 2-D, respectively.
Therefore, while being {\it asymptotically} negligible in determining the rank, for certain configurations, the remainder can appear to govern the rank until higher \(ka\) values are reached. This postpones the eventual asymptotic behavior and creates a false impression of a slower than \(O(N)\) scaling, as  observed in many numerical rank-scaling demonstrations.

The shadow-area based asymptotic rank (or knee index) estimate also guides the physical interpretation of these results. 
It was shown that the predicted knee index partitions, quite accurately, the interacting (or non-silent) domain basis into two subspace bases with different wave characteristics.
The basis vectors up to the knee index,
can be associated with the interacting regions area (or length in 2-D).
For example, for planar domains, they were shown to have relatively uniform spatial intensity distribution and produce/receive directional scanning beams.
As such, they were said to span the ``aperture" subspace of the interaction.

The remaining interacting basis vectors for these cases were shown to be spatially more localized and non-directional, suggesting their association with ``edge diffraction". 
Such locality traits were observed (although not always as pronounced for each individual basis vector) in other configurations of non-flat domains. 
For interactions in the plane in 3-D, 
 the diffraction spectrum itself was shown to have a knee that is geometrically predicted by the shadow length, and a more complex remainder in 3-D than in equivalent 2-D scenarios, alluding to additional diffractive wave constructs.
The interplay of the aperture and diffraction subspaces agrees with the asymptotic and pre-asymptotic behaviors of the rank: 
The diffraction is more dominant for small aspect angles between the source and observer. 
It entirely governs end-fire interactions in quasi-planar configurations. 
This agrees with the favorable compressibility and asymptotic fast solver performance demonstrated for such cases.

Non-asymptotic and precise prediction of the knee index and the rank for arbitrary 3-D configurations requires further investigation. 
In particular, it raises questions on the ability to associate the behavior of the diffraction subspaces with specific wave phenomena and relating them to correction terms in analytical DoF estimators. 
If these can be quantified using some geometric characteristics, this can enable uniform rank prediction, from broadside to end-fire configurations. 
It is also interesting to explore whether savings can be gained by suppression or filtering-out of artificial diffractive behavior that is introduced by the hierarchical partitioning.
 The application of the presented theory to more complicated IE kernels and blocks, including for SIEs for penetrable objects \cite{mavrikakisSurfaceIntegralEquations2025}, in layered  background \cite{liuFastDirectSolution2024} or guiding channels \cite{huangGreensFunctionAnalysis2005}, and with unconventional source distributions \cite{sharshevskyDirectSolutionScattering2020, zvulunGeneralizedSourceIntegral2023,dahanFastDirectSolvers2024, kalhoferFastDirectSolution2025}, opens several additional paths for further research.



\end{document}